\documentclass{article}
\usepackage{amssymb}
\usepackage{multicol}
\usepackage{graphicx}
\usepackage{amsmath, epsfig}
\newcommand{\ignore}[1]{}
\newcommand{\qed}{ \hfill$\square$ }

\begin{document}
\title{Some Ramsey results for the $n$-cube
}
\author{Ron Graham \\
University of California, San Diego
\and Jozsef Solymosi\footnote{This research was supported in part by NSERC and OTKA grants and an Alfred P. Sloan Fellowship. }\\
            University of British Columbia, Vancouver, Canada
            }
            \maketitle
\begin{center}
{\bf Abstract}
\end{center}
In this note we establish a Ramsey-type result for certain subsets of the
$n$-dimensional cube. This can then be applied to obtain reasonable bounds
on various related structures, such as (partial) Hales-Jewett lines for alphabets
of sizes $3$ and $4$, Hilbert cubes in sets of real numbers with small sumsets,
``corners'' in the integer lattice in the plane, and 3-term integer geometric
progressions.
\section{Preliminaries} Let us define
\[
\{0,1\}^n = \{(\epsilon_1, \epsilon_2, \ldots, \epsilon_n), \epsilon_i = 0 \,\mbox{or} \,1, 1 \leq i \leq n\},\quad
D(n):=\{0,1\}^n \times \{0,1\}^n.
\]
We can think of the points of $\{0,1\}^n$ as vertices of an $n$-cube $Q^n$, and $D(n)$ as all the
line segments joining two vertices of $Q^n$. We will ordinarily assume that the two vertices are
distinct. We can represent the points $(X,Y) \in D(n)$ schematically by the diagram
shown in Figure 1.

\begin{figure*}[htbp]
\begin{center}
\includegraphics [scale=.3]{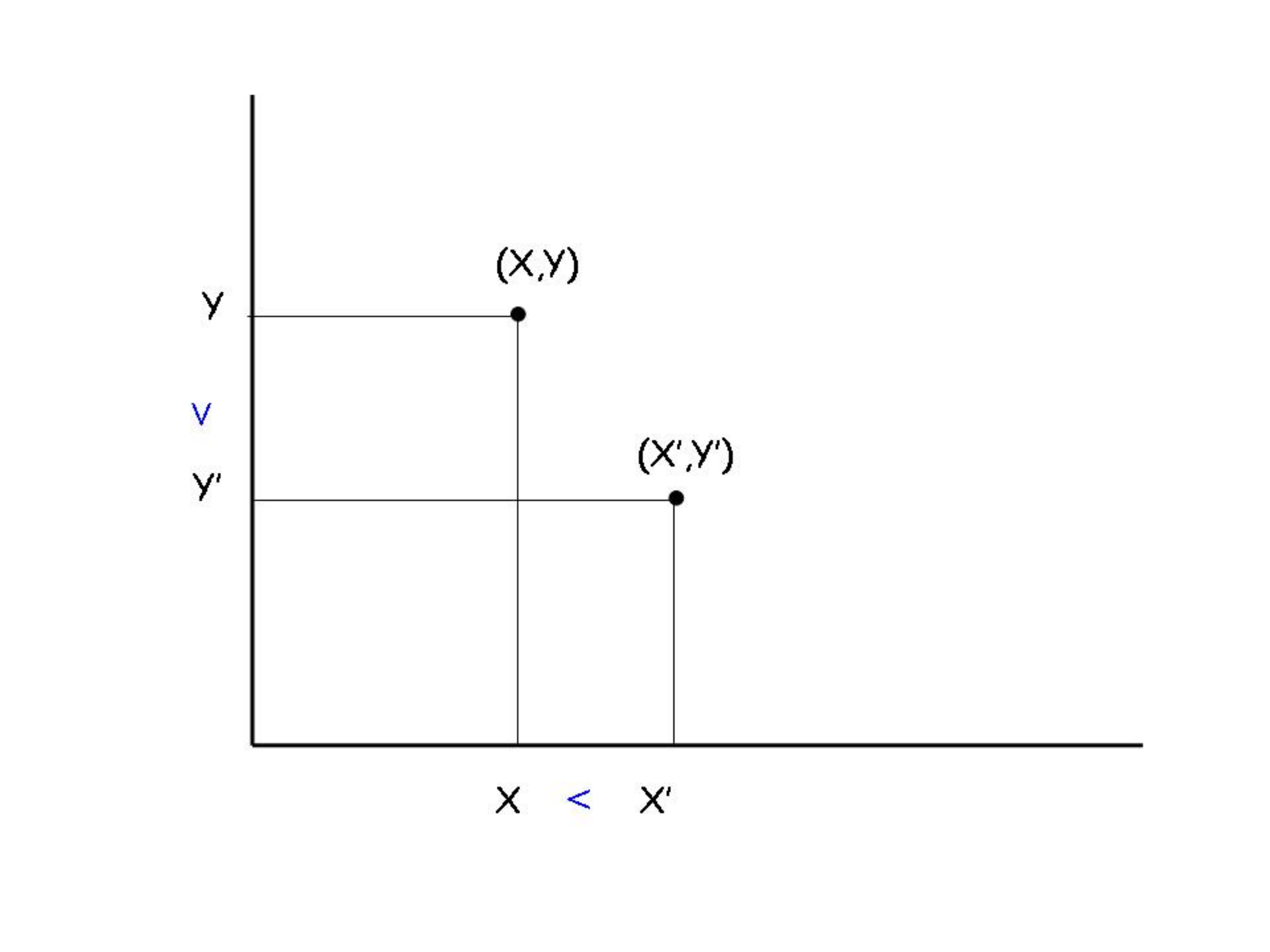}
\end{center}
\caption{Representing points in $D(n)$.}
  \label{fig1}
\end{figure*}

In the diagram, we place $X$ to the left of $X'$ if $w(X) < w(X')$ where $w(Z)$ denotes
the number of $1$'s in the binary $n$-tuple $Z$. Similarly, we place $Y'$ below $Y$ if $w(Y') < w(Y)$. (If $w(X) = w(X')$ or $w(Y) = w(Y')$, then the order doesn't matter).

With $[n]:=\{1, 2,\ldots, n\}$, $I \subseteq [n]$ and $\overline{I} = [n] \backslash I$,
a {\bf line} $L = L(I,C)$ consists all the pairs $((x_1, x_2, \ldots, x_n),(y_1, y_2, \ldots, y_n))$ where $C = (c_j)_{j \in \overline{I}}$ with $y_i = 1-x_i$ if $ i \in I$, and $x_j = y_j = c_j$ if $j \in \overline{I}$.

Thus,
\[
|L(I,C)| = 2^{|I|}.
\]
In this case we say that $I$ has {\it dimension} $|I|$.

\noindent {\bf Fact.} Every point $(X,Y) \in D(n)$ lies on a unique line.\\[.1in]
{\bf Proof}. Just take $I = \{i \in [n]: x_i \neq y_i\}$ and $c_j = x_j = y_j$ for $j \in \overline{I}$.\\

\begin{figure*}[htbp]
\begin{center}
\includegraphics [scale=.3]{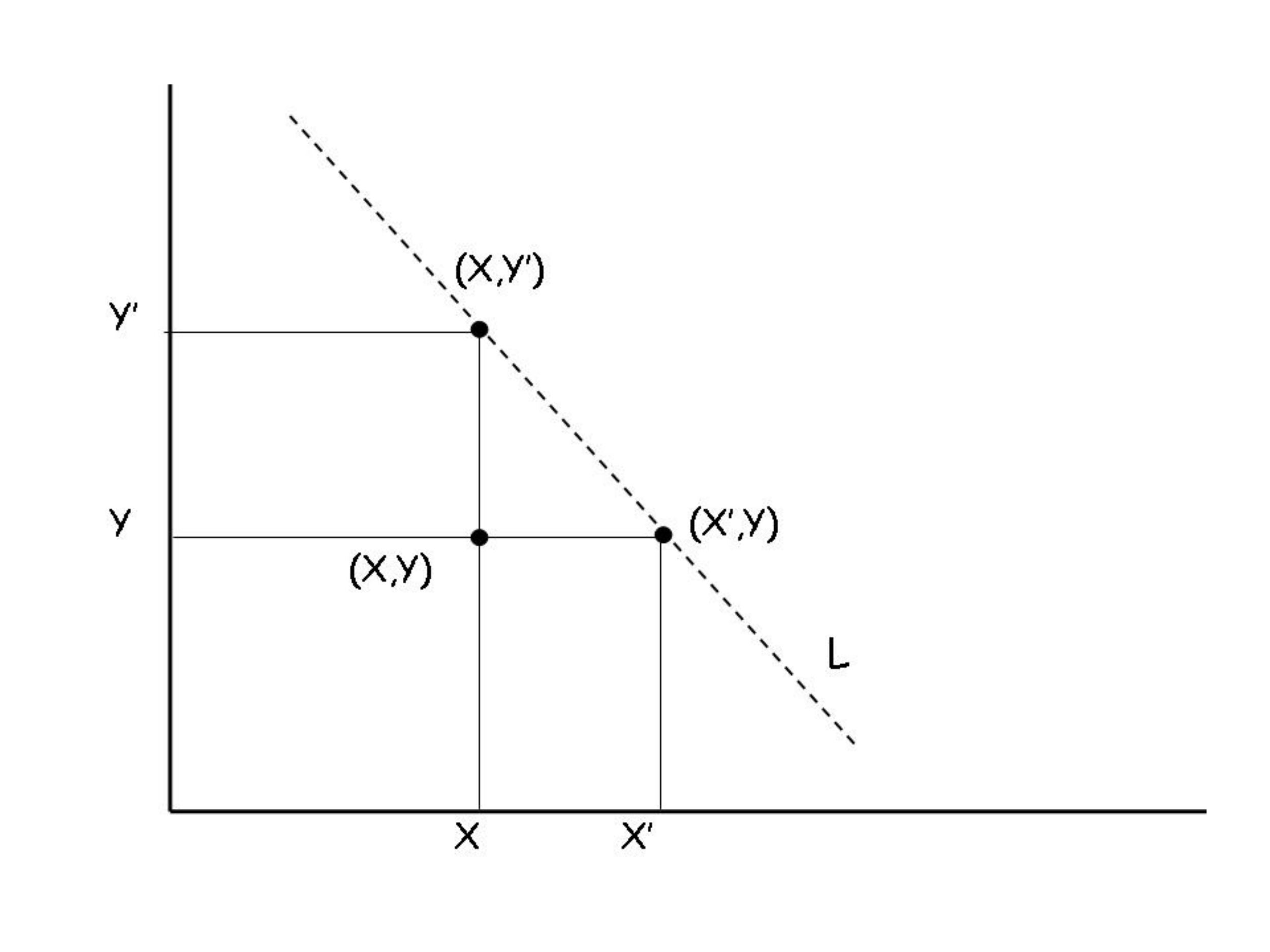}
\end{center}
\caption{A corner in $D(n)$.}
  \label{corner}
\end{figure*}

\noindent By a {\it corner} in $D(n)$, we mean a set of three points of the form
$(X,Y), (X',Y), (X,Y')$ where $(X,Y')$ and $(X',Y)$ are on a common line $L$
(see Figure \ref{corner}).

\begin{figure*}[htbp]
\begin{center}
\includegraphics [scale=.3]{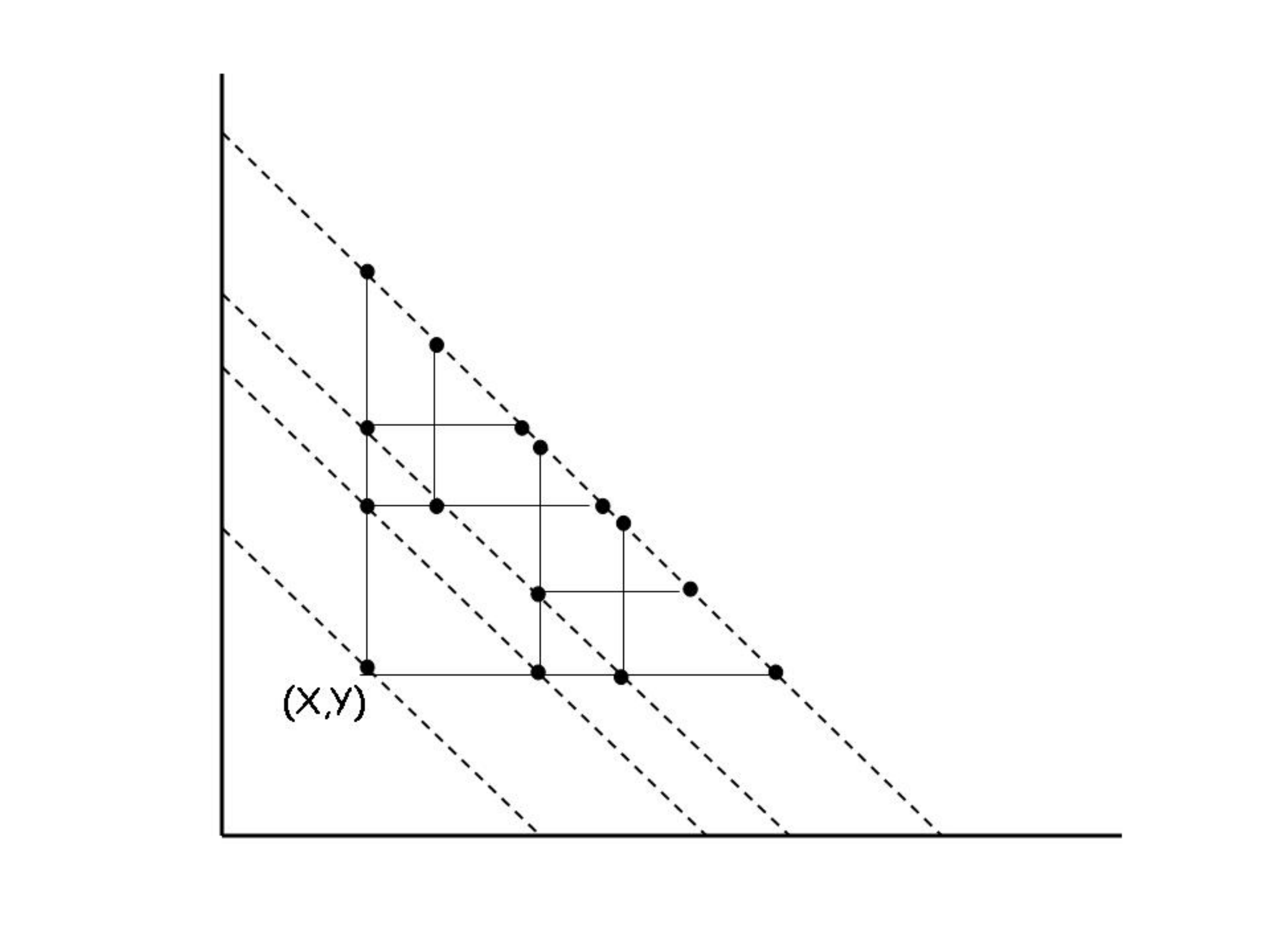}
\end{center}
\caption{A binary tree with 3 levels}
  \label{3_levels}
\end{figure*}

We can think of a corner as a binary tree with
one level and root $(X,Y)$. More generally, a {\it binary tree} $B(m)$ with
$m$ levels and root $(X,Y)$ is defined by joining $(X,Y)$ to the roots of two binary
trees with $m-1$ levels. All of the $2^k$ points at level $k$ are required to be
a common line (see Figure \ref{3_levels}).

\section{The main result.} 
\begin{figure*}[htbp]
\begin{center}
\includegraphics [scale=.3]{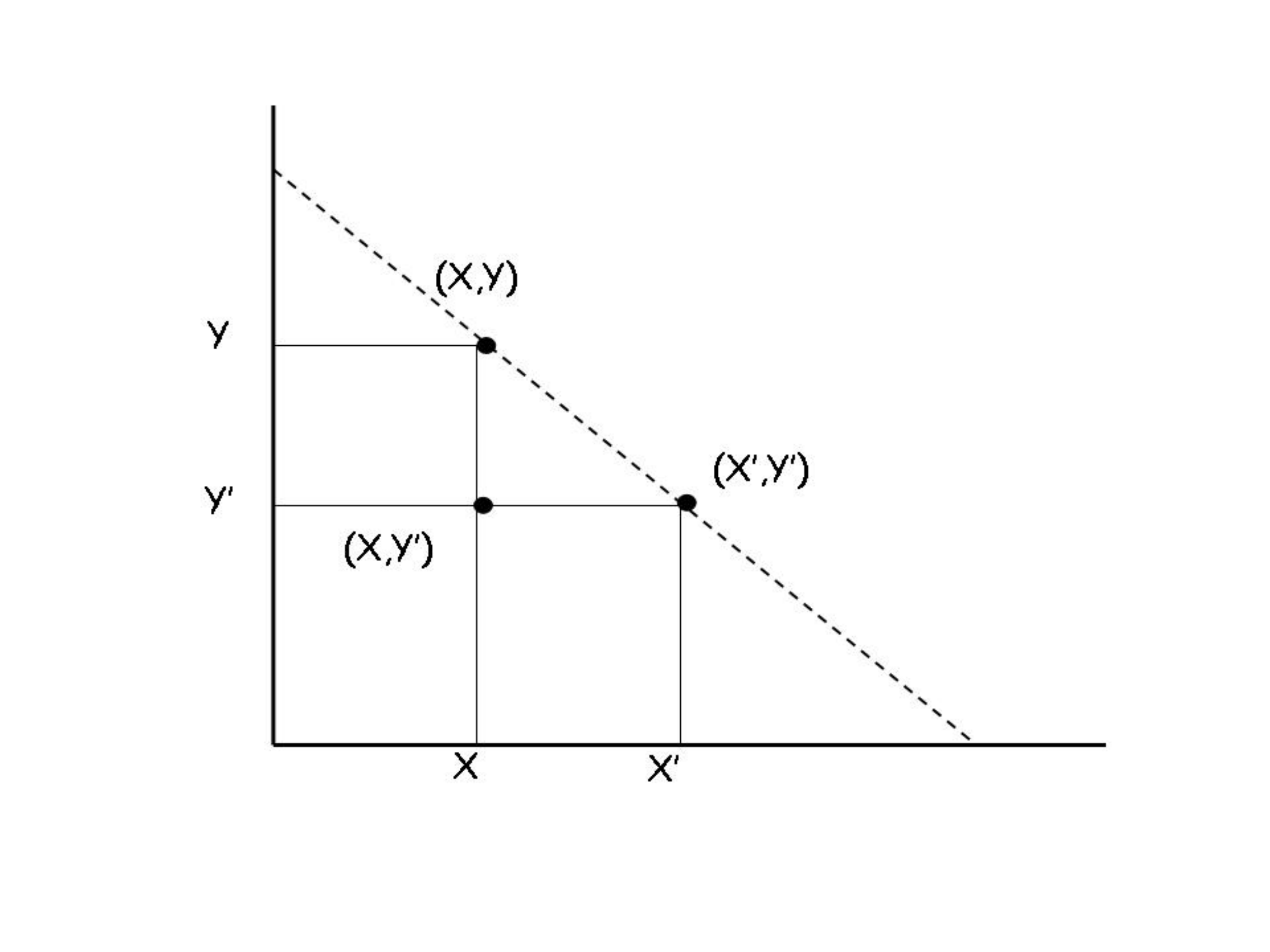}
\end{center}
\caption{A grid point}
  \label{grid}
\end{figure*}

Our first theorem is the following.\\[.2in]
{\bf Theorem 1.} For all $r$ and $m$, there is an $n_0 = n_0(r,m)$ such if $n \geq n_0$
and the points of $D(n)$ are arbitrarily $r$-colored, then there is always a
monochromatic binary tree $B(m)$ with $m$ levels formed. In fact, we can take
$n_0(r,m) = c\,6^{rm}$ for some absolute constant $c$.\\[.1in]
{\bf Proof.} Let $n$ be large (to be specified later) and suppose the points
of $D(n)$ are $r$-colored. Consider the $2^n$ points on the line $L_0 = L([n])$.
Let $S_0 \subseteq L_0 $ be the set of points having the ``most popular'' color $c_0$.
Thus, $|S_0| \geq \frac{2^n}{r}$. Consider the {\it grid} $G_1$
(lower triangular part of a Cartesian
product) defined by:
\[
G_1 = \{(X,Y'): (X,Y) \in S_0, (X', Y') \in S_0 \, \mbox{with} \, X < X'\}.
\]
(See Figure \ref{grid}).

\noindent Thus,
\[
|G_1| \geq \left( S_0 \atop 2 \right) > \frac{1}{4} |S_0|^2 \geq \frac{1}{4r^2} \cdot 4^n:=\alpha_1 4^n.
\]

Let us call a line $L$ of dimension $t$ {\it small} if $t < \frac{n}{3}$ and
{\it deficient} if \\$|L \cap G_1| \leq (\frac{\alpha_1}{4}) 2^t$.\\
Thus, the total number of points on small or deficient lines is at most
\begin{eqnarray*}
 && \sum_{t<\frac{n}{3}} 2^t \left( n \atop t \right) 2^{n-t} + \sum_{t \geq \frac{n}{3}} \frac{\alpha_1}{4} 2^t \left( n \atop t \right) 2^{n-t}\\
&=& \frac{\alpha_1}{4} \sum_t 2^n \left( n \atop t \right) + (1 - \frac{\alpha_1}{4})
\sum_{t < \frac{n}{3}} 2^n \left( n \atop t \right)\\
& \leq & (\frac{\alpha_1}{4}) 4^n + (1 - \frac{\alpha_1}{4})(3.8^n)\\
&&\quad \quad \quad \quad (\mbox{since} \sum_{t<\frac{n}{3}} \left( n \atop t \right) <  1.9^n \,\mbox{follows easily by induction})\\[-.2in]
& \leq & (\frac{\alpha_1}{2}) 4^n
\end{eqnarray*}
provided $\alpha_1 \geq 2 \cdot (.95)^n$.\\

\begin{figure*}[htbp]
\begin{center}
\includegraphics [scale=.3]{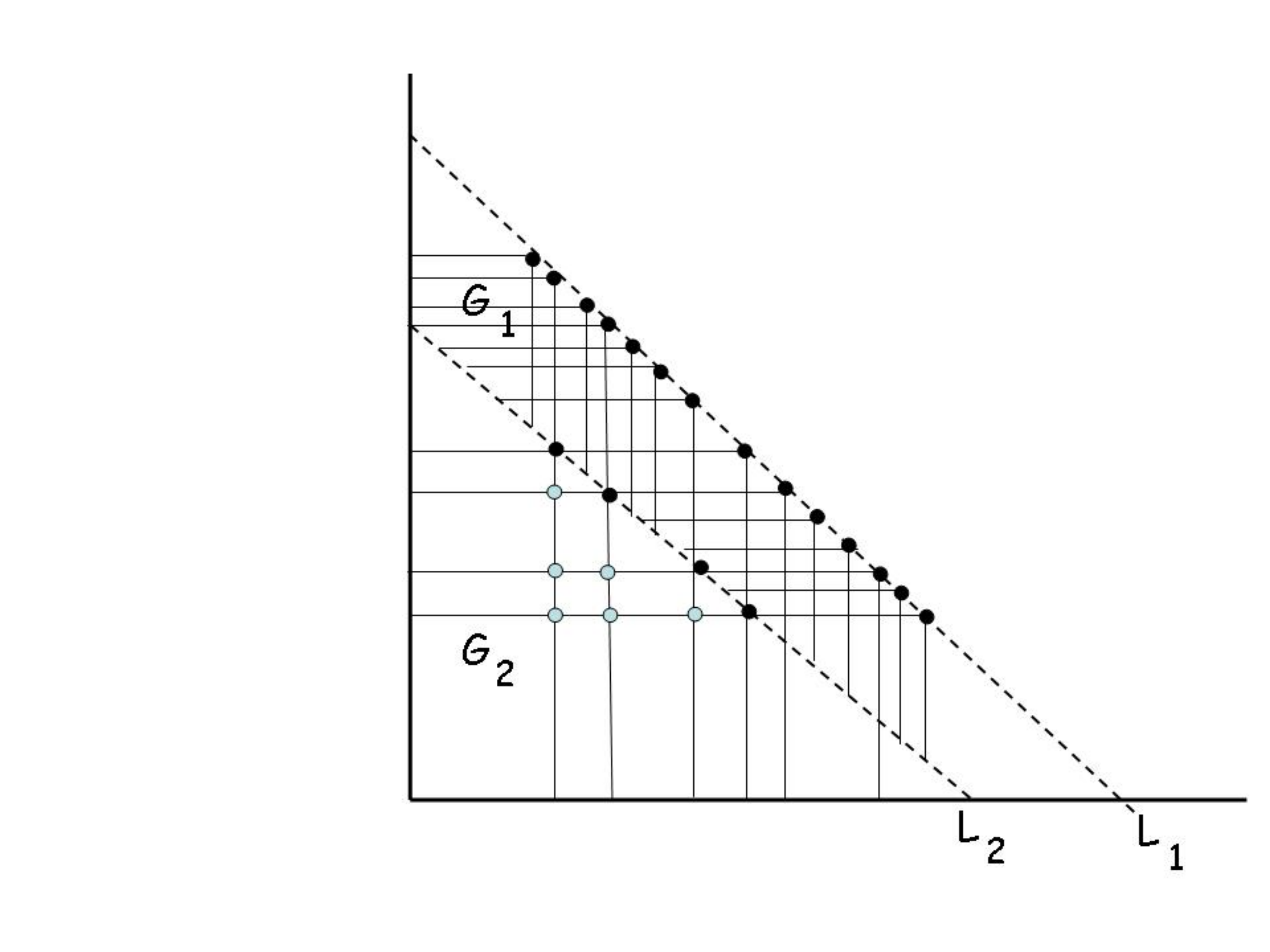}
\end{center}
\caption{$G_2 \subset G_1$}
  \label{G1_in_G2}
\end{figure*}

Thus, if we discard these points, we still have at least $(\frac{\alpha_1}{2})4^n$ points
remaining in $G_1$, and all these points are on ``good'' lines, i.e., not small
and not deficient.
\noindent Let $L_1$ be such a good line, say of dimension $|I_1| = n_1 \geq \frac{n}{3}$. Let $S_1$ denote the set of points of $L_1 \cap G_1$ with the
most popular color $c_1$. Therefore
\[
|S_1| \geq (\frac{\alpha_1}{4r}) 2^{n_1}.
\]
Observe that $G_2 \subset G_1$ (see Figure \ref{G1_in_G2}).

\noindent Now let $G_2$ denote the ``grid'' formed by $S_1$, i.e.,
\[
G_2 = \{(X,Y'): (X,Y \in S_1, (X',Y') \in S_1, \, \mbox{with} \, X < X'\}.
\]
Therefore, we have
\[
|G_2| \geq \left(|S_1| \atop 2 \right) \geq (\frac{\alpha_1}{8r})^2\,4^{n_1} :=\alpha_2 4^{n_1}.
\]

\noindent As before, let us classify a line $L$ of dimension $t$ as {\it small}
if $t < \frac{n_1}{3}$, and as {\it deficient} if $|L \cap G_2| \leq (\frac{\alpha_2}{4}) 2^t.$

\noindent A similar calculation as before shows that if we remove from $G_2$ all
the points on small or deficient lines, then at least $(\frac{\alpha_2}{2})4^{n_1}$
points will remain in $G_2$, provided $\alpha_2 \geq 2 \cdot (.95)^{n_1}.$

\noindent Let $S_2 \subseteq L_2 \cap G_2$ have the most popular color $c_2$, so that
\[
|S_2| \geq (\frac{\alpha_2}{4r})2^{n_2}.
\]
Then, with $G_3$ defined to the the ``grid'' formed by $S_2$, we have
$|G_3| \geq (\frac{\alpha_2}{8r})^2 4^{n_2}$, and so on. Note that
$G_3 \subset G_2 \subset G_1$.

\noindent We continue this process for $rm$ steps.\\
In general, we define
\[
\alpha_{i+1} = (\frac{\alpha_i}{8r})^2, i \leq i \leq rm -1
\]
with $\alpha_1 = \frac{1}{4r^2}$.  By construction, we have $n_{i+1} \geq \frac{n_i}{3}$
for all $i$. In addition, we will need to have $\alpha_i \geq 2 \cdot (.95)^{n_i}$ for
all $i$ for the argument to be valid. In particular, this implies that in general
\[
\alpha_k = \frac{1}{2^{2^{k+2}- 6} r^{2^{k+1} - 2}}.
\]

\noindent It is now straightforward to check that all the required inequalities
are satisfied by choosing $n \geq n_0(r,m) = c \cdot 6^{rm}$ for a suitable
absolute constant $c$.

\noindent Hence, there must be $m$ indices $i_1 < i_2 < \ldots < i_m$ such
that all the sets $S_{i_k}$ have the same color.

\noindent These $m$ sets $S_{i_k}$ contain the desired monochromatic binary tree $B(m)$.

\section{Some interpretations}
\subsection{Self-crossing paths} As we stated at the beginning, we can think
of $D(n)$ as the set of all the diagonals of the $n$-cube $Q^n$. Let us
call a pair $\{x,\bar x\} = \{(x_1, \ldots, x_n),(\bar x_1, \ldots, \bar x_n)\}
$ a {\bf main} diagonal of $Q^n$ where $\bar x_i = 1 - x_i$.

\noindent An affine $k$-subcube of $Q^n$ is defined to be a subset of $2^k$
points of the form
$\{(y_1, \ldots, y_n) : y_i = 0$ or $1$ if and only if $i \in I \}$ for
some $k$-subset $I \subseteq [n] = \{1,2, \ldots, n\}$.

\noindent We will say that three connected diagonals of the form
$\{x,y\}, \{y,z\}, \{z,w\}$ form a {\it self-crossing path}, denoted by $\ltimes$,
if $\{x,y\}$ and $\{z,w\}$ are both main diagonals of the same subcube.

\noindent {\bf Corollary 1}. In any $r$-coloring of the edges in $D(n)$, there is
always a monochromatic self-crossing path $\ltimes$, provided $n > c \cdot 6^r$
(where $c$ is a suitable absolute constant).

\noindent The same argument works for any subgraph $G$ of $D(n)$, provided
that $G$ has enough edges and for any pair of crossing main diagonals, $G$ has all
the edges between the pair's endpoints.

\subsection{Corners}
The preceding techniques can be used to prove the following.\\
{\bf Theorem 3.} For every $r$, there exists $\delta = \delta(r)$ and
$n_0 = n_0(r)$ with the following property:\\[.2in]
If $A$ and $B$ are sets of real numbers with $|A| = |B| = n \geq n_0$ and
$|A + B| \leq n^{1+ \delta}$, then any $r$-coloring of $A \times B$ contains
a monochromatic ``corner'', i.e., a set of $3$ points of the form
$(a,b), (a+d,b), (a,b+d)$ for some positive number $d.$
In fact, the argument shows that we can choose $\delta = \frac{1}{{2^{r+1}}}$.


The calculation goes as follows; The Cartesian product $A\times B$ can be
covered by $n^{1+\delta}$ lines of slope -1. Choose the line with the most
points from $A\times B,$ denoted by $L_0.$ There are at least $n^2/n^{1+\delta}$
points in $L_0\cap A\times B.$ Choose the set of points $S_1$ with the most popular color
in $L_0\cap A\times B.$  ($|S_1|\geq n^{1-\delta}/r$) As before, consider the grid $G_1$
defined by $S_1,$ and choose the slope -1 line, $L_2,$
which has the largest intersection with $G_1.$ Choose the set of points, $S_2,$ having the most
popular color  and repeat the process with $G_2,$ the grid defined by $S_2.$   We can't
have more than $r$ iterations without having a monochromatic corner. Solving the simple recurrence
$a_{n+1}=2a_n+1$ in the exponent, one can see that after $r$ steps
the size of $S_r$ is at least $c_rn^{1-\delta(2^{r+1}-1)}.$ If this quantity is at least 2, then we have
at least one more step and the monochromatic corner is unavoidable. The inequality
\[
c_rn^{1-\delta(2^{r+1}-1)} \geq 2
\]
can be rearranged into
\[
n^{1-\delta{2^{r+1}}}\geq {2\over{c_rn^{\delta}}} .
\]
From this we see that choosing $\delta=2^{-r-1}$ guarantees that for large enough $n$
the inequality above is true, proving our statement.


By iterating these techniques, one can show that the same hypotheses on $|A|$
and $|B|$ (with appropriate $\delta = \delta(r,m)$ and $n_0 = n_0(r,m)$, imply that if
$A \times B$ is $r$-colored then each set contains a monochromatic translate
of a large ``Hilbert cube'', i.e., a set of the form
\[
H_m(a,a_1, \ldots, a_m) = \{a + \sum_{1 \leq i \leq m} \epsilon_i a_i \} \subset A,
\]
\[
H_m(b,a_1, \ldots, a_m) = \{b + \sum_{1 \leq i \leq m} \epsilon_i a_i \} \subset B
\]
where $\epsilon_i = 0$ or $1, \, 1 \leq i \leq m$.

\subsection{Partial Hales-Jewett lines}
{\bf Corollary 2.} For every $r$ there is an $n = n_0(r)\leq c6^r,$ with the following
property. For every $r$-coloring of $\{0,1,2,3\}^n$ with $n > n_0$, there is
always a monochromatic set of $3$ points of the form:
\begin{eqnarray*}
(\ldots,a,\ldots,0,\ldots,b,\ldots,3,\ldots,0,\ldots,c,\ldots,3,\ldots,d,\ldots)\\
(\ldots,a,\ldots,1,\ldots,b,\ldots,2,\ldots,1,\ldots,c,\ldots,2,\ldots,d,\ldots)\\
(\ldots,a,\ldots,2,\ldots,b,\ldots,1,\ldots,2,\ldots,c,\ldots,1,\ldots,d,\ldots)
\end{eqnarray*}
In other words, every column is either {\it constant}, {\it increasing} from $0$,
or {\it decreasing} from 3.\\
{\bf Proof.}  To each point $(x_1, x_2, \ldots, x_n)$ in $\{0,1,2,3\}^n$, we
associate the point $\big( (a_1, a_2, \ldots, a_n), (b_1, b_2, \ldots, b_n) \big)$
in $\{0,1\}^n \times \{0,1\}^n$ by the following rule:
\[
\begin{array}{cccc}
x_k & \leftrightarrow & a_k & b_k \\ \hline
0 & \vline & 0 & 0 \\
1 & \vline & 0 & 1 \\
2 & \vline & 1 & 0 \\
3 & \vline & 1 & 1
\end{array}
\]
Then it not hard to verify that a monochromatic corner in $D(n)$ corresponds to
a monochromatic set of $3$ points as described above, a structure which we might
call a partial Hales-Jewett line.

\medskip

\noindent
{\bf Corollary 3.} For every $r$ there is an $n = n_0(r)\leq c6^r,$ with the following
property. For every $r$-coloring of $\{0,1,2\}^n$ with $n > n_0$, there is
always a monochromatic set of $3$ points of the form:
\begin{eqnarray*}
(\ldots,a,\ldots,0,\ldots,b,\ldots,0,\ldots,0,\ldots,c,\ldots,0,\ldots,d,\ldots)\\
(\ldots,a,\ldots,1,\ldots,b,\ldots,2,\ldots,1,\ldots,c,\ldots,2,\ldots,d,\ldots)\\
(\ldots,a,\ldots,2,\ldots,b,\ldots,1,\ldots,2,\ldots,c,\ldots,1,\ldots,d,\ldots)
\end{eqnarray*}
{\bf Proof.}  Map the points $(a_1, a_2, \ldots, a_n) \in \{0,1,2,3\}$ to \\points
$(b_1, b_2, \ldots, b_n) \in \{0,1,2\}^n$ by:\\[-.3in]
\begin{center}
$a_i = 0$ or $3 \Rightarrow b_i = 0$, $a_i = 1 \Rightarrow b_i = 1$, $a_i = 2.
 \Rightarrow b_i = 2$
\end{center}
The theorem now follows by applying Corollary 2. \qed

\subsection{3-term geometric progressions.}
The simplest non-trivial case of  van der Waerden's theorem \cite{vW} states that
for any natural number $r$, there is a number $W(r)$ such that for any
$r$-coloring of the first $W(r)$ natural numbers there is a monochromatic
three-term arithmetic progression. Finding the exact value of $W(r)$ for
large $r$-s is a hopelessly difficult task. The best upper bound follows
from a recent result of Bourgain \cite{Bo}; \[W(r)\leq ce^{r^{3/2}}.\]
One can ask the similar problem for geometric progressions;
What is the maximum number of colors, denoted by $r(N),$
that for any $r(N)$-coloring of the first $N$ natural numbers
there is a monochromatic geometric progression. Applying
Bourgain's bound to the exponents of the geometric progression
$ \{2^i\}_{i=0}^{\infty},$ shows that $r(N)\geq c\log\log{N}.$
Using our method we can obtain the same bound, without applying
Bourgain's deep result.

\medskip

\noindent
Observe that if we associate the point
$(a_1, a_2,\ldots,a_k, \ldots, a_n)$ with the integer $\prod_k p_k^{a_k}$, where
 $p_i$ denotes the $i^{th}$ prime, then the points
\begin{eqnarray*}
(\ldots,a,\ldots,0,\ldots,b,\ldots,3,\ldots,0,\ldots,c,\ldots,3,\ldots,d,\ldots)\\
(\ldots,a,\ldots,1,\ldots,b,\ldots,2,\ldots,1,\ldots,c,\ldots,2,\ldots,d,\ldots)\\
(\ldots,a,\ldots,2,\ldots,b,\ldots,1,\ldots,2,\ldots,c,\ldots,1,\ldots,d,\ldots)
\end{eqnarray*}
correspond to a $3$-term geometric progression. Our bound from
Corollary 2 with an estimate for the product of the first $n$ primes
imply that $r(N)\geq c\log\log{N}.$

\section{Concluding remarks}
It would be interesting to know if we can ``complete the square'' for some of
these results. For example, one can use these methods to show that if the points
of $[N] \times [N]$ are colored with at most $c \log \log N$ colors, then there
is always a monochromatic ``corner'' formed, i.e., $3$ points $(a,b), (a',b), (a,b')$
with $a' + b = a + b'$. By projection, this gives a $3$-term arithmetic progression
(see \cite{gs}).

\noindent Is it the case that with these bounds (or even better ones), we can
guarantee the $4^{th}$ point $(a',b')$ to be monochromatic as well?

\noindent Similarly, if the diagonals of an $n$-cube are $r$-colored
with $r < c \log \log n$,
is it true that a monochromatic $\bowtie$ must be formed, i.e., a self-crossing $4$-cycle
(which is a self-crossing path with one more edge added)?

\noindent Let $\boxtimes$ denote the structure consisting of the set of $6$ edges spanned
by $4$ {\bf coplanar} vertices of an $n$-cube. In this case, the occurrence
of a monochromatic $\boxtimes$ is guaranteed once $n \geq N_0$, where $N_0$
is a {\bf very} large (but well defined) integer, sometimes referred to as Graham's number
(see \cite{wiki}). The best lower bound currently available for $N_0$ is $11$
(due to G. Exoo \cite{exoo}).

One can also ask for estimates for the density analogs for the preceding results.
For example, Shkredov has shown the following:\\
{\bf Theorem} (\cite{shkredov}). Let $\delta > 0$ and $N \ll \exp \exp(\delta^{-c})$,
where $c > 0$ is an absolute constant. Let $A$ by a subset of $\{1, 2,\ldots, N\}^2$
of cardinality at least $\delta N^2$. Then $A$ contains a corner.

\noindent It would be interesting to know if the same hypothesis implies that
$A$ contains the $4^{th}$ point of the corner, for example.

\end{document}